\documentclass[12pt]{article}
\usepackage{amsmath,amsfonts,amssymb, color}
\date{}
\newtheorem{Theorem}{Theorem}[section]

\newtheorem{Proposition}[Theorem]{Proposition}
\newtheorem{Lemma}[Theorem]{Lemma}
\newtheorem{Corollary}[Theorem]{Corollary}
\newtheorem{Remark}[Theorem]{Remark}

\newtheorem{Hypothesis}[Theorem]{Hypothesis}
\makeatletter
\newif\ifmsbmloaded@

\@addtoreset{equation}{section}

\makeatother

\newcommand{\one}{1\!\!\!\;\mathrm{l}}

\def\R{\mathbb R}
\def\N{\mathbb N}

\def\P{\mathbb P}

\def\ds{\displaystyle}

\title{\bf Uniqueness for solutions of Fokker--Planck equations on infinite dimensional spaces}
\author{Vladimir Bogachev\thanks{Supported in part by
the RFBR project
 07-01-00536,
the Rus\-sian--Japanese Grant 08-01-91205-JF, the
Rus\-sian--Ukrainian Grant 08-01-90431,
 SFB 701 at the University of Bielefeld.},\\
Department of Mechanics and Mathematics,\\
          Moscow State University, 119991 Moscow, Russia,\\\\
Giuseppe Da Prato,\\
 Scuola Normale Superiore
di Pisa, Italy\\
and\\
 \\
Michael R\"ockner \thanks{Supported by the DFG through SFB-701 and IRTG 1132 as well as the
BIBOS-Research Center. } \\
Faculty of Mathematics, University of Bielefeld, Germany\\ and\\
Department of   Statistics,\\ Purdue University, W. Lafayette, 47906, IN,   U. S. A.
}

\begin{document}

\maketitle
\begin{abstract}
We develop a general technique to prove uniqueness of solutions for Fokker--Planck equations on infinite dimensional spaces. We illustrate this method by implementing it for Fokker--Planck equations in Hilbert spaces with
Kolmogorov operators with irregular coefficients and both non-degenerate or degenerate second order part.

 \end{abstract}

\noindent {\bf 2000 Mathematics Subject Classification AMS}: 60H15, 60J35, 60J60, 47D07


\bigskip

\maketitle

\section{Introduction}
 Fokker--Planck and transport equations with
irregular coefficients in finite dimensions  have been studied intensively in recent years (see e.g.
\cite{A04}, \cite{A05}, \cite{BDPR04}, \cite{BDPR08}, \cite{BDPRS07}, \cite{BRS09}, \cite{Fi08}, \cite{LeLi08},  \cite{L07}, \cite{M07}
  and the references therein,  and also the fundamental paper \cite{DPL89}).
More recently   transport  and Fokker--Planck equations  have also been studied in infinite dimensions
(see, e.g., \cite{0}, \cite{BDPRS09} and  \cite{BDPR08b}, \cite{BDPR08c}, \cite{BDPR09}, \cite{BDPR10}
respectively)

In this paper we further develop the method from \cite{BDPR09}, \cite{BDPR10} to prove uniqueness of solutions of Fokker--Planck equations in infinite dimensions. Our main aim here is to give an independent general   presentation of the single steps and to implement this method under considerably weakened assumptions
 on the coefficients. Though this method to prove uniqueness is more universal and  can be applied to more general Kolmogorov operators, as e.g. those of nonlocal (i.e. pseudo-differential) type, thus allowing jumps for the corresponding stochastic dynamics, we here confine ourselves to the case where, at least  on a heuristic level, there is an underlying  stochastic differential equation in the background.
 More precisely our framework is as follows:\medskip

Let $H$ be a separable real Hilbert space with
inner product $\langle\cdot,\cdot   \rangle$ and corresponding norm $|\cdot|$.  $L(H)$ denotes the set of all bounded linear  operators  on $H$ with its usual norm $\|\cdot\|$, $ \mathcal B(H)$ its Borel $\sigma$-algebra, $   \mathcal B_b(H)$ the set of all bounded $ \mathcal B(H)$-measurable functions
from $H$ to $\R$  and $ \mathcal P(H)$ the set of all probability measures on $H$, more precisely on $(H, \mathcal B(H))$.

Consider the following type of non-autonomous stochastic differential equations on $H$ and time interval $[0,T]$:
\begin{equation}
\label{e1.1}
\left\{\begin{array}{l}
dX(t)=(AX(t)+F(t,X(t)))dt+\sqrt{C}dW(t),\\
\\
X(s)=x\in H,\;t\ge s.
\end{array}\right.
\end{equation}
Here  $W(t),\;t\ge 0,$ is a cylindrical Wiener process on $H$ defined on a stochastic basis $(\Omega, \mathcal F, (\mathcal F_t)_{t\ge 0},\P)$, $C$ is a symmetric positive operator in $L(H)$,
$D(F)\in\mathcal B([0,T]\times H),$
$F\colon D(F)\subset [0,T]\times H\to H$  is a measurable map, and
$A\colon\, D(A)\subset H\to H$ is the infinitesimal generator of a $C_0$-semigroup $e^{tA},\;t\ge 0,$ in $H.$

Without further regularity assumptions on $F$ it is, of course, not at all clear whether \eqref{e1.1} has a solution in the strong or even in the weak sense.
If, however,
there is a weak solution to \eqref{e1.1}, then it is a well known consequence of It\^o's formula that its transition probabilities $p_{s,t}(x,dy),\;x\in H,\;s\le t$, solve the Fokker--Planck equation determined by  the associated Kolmogorov operator, see e.g. \cite{DPZ1}.
But as shown in   our earlier papers \cite{BDPR08b}, \cite{BDPR09}, \cite{BDPR10} one can describe very general conditions on $F$ above for which one can solve the Fokker--Planck equation directly for Dirac initial conditions and thus to obtain the transition functions $p_{s,t},\;s\le t$, corresponding to \eqref{e1.1} though one might not have a solution to it.

The general motivation to study Fokker--Planck equations instead of Kolmogorov equations, as done in some of our former  papers (see e.g. \cite{BRS00}, \cite{BRS02}, \cite{DP2004}, \cite{6}, \cite{DPT01}  and the references therein)   is that the latter are equations for functions, whereas the first are equations for measures for which one has e.g. much better compactness criteria in our infinite dimensional situation. So, there is a good chance to obtain very general existence results. Uniqueness, however, is considerably harder to prove and this is in the centre of considerations in this paper.
\medskip

Before we write down the Fokker--Planck equation precisely we
recall that the Kolmogorov operator $L_0$ corresponding to \eqref{e1.1} reads as follows:
\begin{multline}\label{e1.1'}
L_0u(t,x)=D_tu(t,x)+\frac12\;\mbox{\rm Tr}\;[CD_x^2u(t,x)]\\+\langle x,A^*D_xu(t,x)
 \rangle+\langle F(t,x),D_xu(t,x)  \rangle,\quad x\in H,\;t\in[0,T],
\end{multline}
where $D_t$ denotes the derivative in time and $D_x,D^2_x$ denote the first and second order
Fr\'echet derivatives in space, i.e. in $x\in H$,  respectively.
The operator $L_0$ is defined  on the space $D(L_0):={\mathcal E}_{A}([0,T]\times H)$,   the linear span of all
 real   parts of    functions $u_{\phi,h}$ of the form
 \begin{equation}
\label{e1.3}
u_{\phi,h}(t,x)=\phi(t)e^{i\langle x,h(t) \rangle},\quad t\in [0,T], \;x\in H,
\end{equation}
where $\phi\in C^1([0,T]),\;\phi(T)=0$, $h\in C^1([0,T];D(A^*))$  and $A^*$ denotes the adjoint of $A$.

For a fixed initial time $s\in[0,T]$ the  Fokker--Planck equation   is an equation for measures  $\mu(dt,dx)$ on $[s,T]\times H$ of the type
\begin{equation}\label{e1.2'}
\mu(dt,dx)=\mu_t(dx)dt,
\end{equation}
with $\mu_t\in \mathcal P(H)$ for all $t\in[s,T],$ and $t\mapsto \mu_t(A)$ measurable on $[s,T]$ for all $A\in \mathcal B(H)$,
i.e., $\mu_t(dx),\;t\in[s,T],$ is a probability kernel from
$([s,T],\mathcal B([s,T])$ to $(H,\mathcal B  (H))$. Then the equation
for an initial condition $\zeta\in  \mathcal P(H)$ reads as follows:
$\forall\;u\in D(L_0)$ one has
\begin{multline}
\label{e1.5}
\int_Hu(t,y)\mu_t(dy)=\int_Hu(s,y)\zeta(dy)+\int_s^tds' \int_HL_0u(s',y)\mu_{s'}(dy)
,\\\quad \mbox{\rm for   $dt$-a.e.}\;t\in [s,T],
\end{multline}
where the   $dt$-zero set may depend on $u$.
When writing \eqref{e1.5}   (or \eqref{(1.8)} below) we always implicitly assume that
\begin{equation}
\label{e1.5'}
\int_{[0,T]\times H}(|\langle y,A^*h(t) \rangle|+|F(t,y)|)\mu(dt,dy)<\infty
\end{equation}
for all $h\in C^1([0,T];D(A^*))$ with $|F(t,y)|:=+\infty$ if $(t,y)\notin D(F)$, so that all involved integrals exist in the usual sense.
\begin{Remark}
\label{r1.1}
\em
(i) Considering $D(L_0)$ as test functions and dualizing  it is easy to see that   \eqref{e1.5}  turns into the more familiar form of the Fokker--Planck equation
 \begin{equation}
 \label{(1.7)}
\frac{\partial}{\partial t}\;\mu_t=-L_0^*\mu_t,\quad \mu_s=\zeta.
\end{equation}

\noindent(ii) Setting $t=T$ and recalling that $u(T,\cdot)\equiv 0$ for all $u\in D(L_0)$ we see that
(under assumption \eqref{e1.5'}) equation \eqref{e1.5} is obviously   equivalent to
\begin{equation}\label{(1.8)}
\int_{[s,T]\times H}L_0 u(s',y)\mu(ds',dy)=-\int_Hu(s,y)\zeta(dy),\quad\forall\;u\in D(L_0).
\end{equation}

\end{Remark}
Solving \eqref{e1.5} (if this is possible) with $\zeta=\delta_x$ (:=Dirac measure in $x\in H$) for $x\in H$ and $s\in [0,T)$ and expressing the dependence on $x,s$ in the notation, we obtain probability measures $p_{s,t}(x,dy),\;t\in[s,T],$
such that the measure $p_{s,t}(x,dy)dt$ on $[s,T]\times H$ is a solution of \eqref{e1.5}. It was proved in detail in Section 3 of \cite{BDPR10} that if we have
uniqueness for \eqref{e1.5} and  ``sufficient continuity'' of the functions
 $t\mapsto p_{s,t}(x,dy)$, then these measures satisfy the Chapman--Kolmogorov equations, i.e. for $0\le r<s<t\le T$ and $x\in H$ (or in a properly chosen subset thereof)
\begin{equation}
\label{e1.9}
\int_Hp_{s,t}(x',dy)p_{r,s}(x,dx')=p_{r,t}(x,dy),
\end{equation}
where the left hand side is a measure defined for $A\in \mathcal B(H)$ as
$$
\int_{H\times H}\one_A(y)p_{s,t}(x',dy)p_{r,s}(x,dx').
$$
 \medskip

 In all of this paper we shall concentrate on conditions on the coefficients $A$, $C$ and $F$  in \eqref{e1.1} under which we can prove uniqueness, not caring about existence at all, since the last was studied in detail in \cite{BDPR08b}, \cite{BDPR09}  and\cite{BDPR10}. Unlike in the previous work we include both cases with Tr $C=+\infty$ and
Tr $C<+\infty$.\medskip

The organization of the paper is as follows:

In Section 2 we explain the general argument, namely that ``the dense range condition'' (cf. \eqref{e2.1} below) implies uniqueness of solutions to \eqref{e1.3}.

In the subsequent sections we show how to check the ``the dense range condition''.
To this end in Section 3 we recall some known regularity results for the time dependent Ornstein--Uhlenbeck operator on Hilbert spaces from \cite{BDPR08c} and some of its consequences to be used below. In Section 4 we show that ``the dense range condition'' holds and hence that \eqref{e1.3} has at most one solution in the case $C^{-1}\in L(H)$. This can be done just  under an $L^2$-integrability condition on $F$. Section 5 is devoted to possibly degenerate cases where not necessarily $C$ is invertible, including the deterministic case where $C=0$. Section 6 contains applications.\medskip

Finally, we would like to mention that even if \eqref{e1.1} has a unique solution, it is not clear at all why the Fokker--Planck equation \eqref{e1.5} has a unique solution. For instance, there could be solutions $p_{s,t}(x,dy)$ to \eqref{e1.5} for $\zeta=\delta_x$ for every $x\in H$ satisfying the Chapman--Kolmogorov equation \eqref{e1.9}, for which there exists no process with continuous or cadlag paths so that $p_{s,t}(x,dy),\;x\in H,\;s\le t,$ are its transition probabilities.

\section{The general argument}

Fix $\zeta\in\mathcal P(H)$, $s\in[0,T]$ and define
$\mathcal M_{s,\zeta}$ the set of all finite nonnegative measures $\mu(dt\,dx)$ on $[s,T]\times H$ satisfying \eqref{e1.2'}, \eqref{e1.5} and \eqref{e1.5'}.

Then we have the following general result:
\begin{Theorem}
\label{t2.1}
Let $\mathcal K\subset\mathcal M_{s,\zeta}$ be a convex subset such that the following
``dense range condition'' is satisfied
\begin{equation}
\label{e2.1}
L_0(D(L_0))\;\mbox{\it is dense in }\;L^1([0,T]\times H,\mu),
\end{equation}
for all $\mu\in \mathcal K.$ Then $\mathcal K$ contains at most one element.
\end{Theorem}
{\bf Proof}. Let $\mu^{(i)}(dt\,dx)=\mu_t^{(i)}(dx)dt\in \mathcal K,\;i=1,2.$ Then
\begin{equation}
\label{e2.2}
\mu_t (dx)dt:=\frac12\;\mu_t^{(1)}(dx)dt+\frac12\;\mu_t^{(2)}(dx)dt\in \mathcal K
\end{equation}
and any $\mu(dt\,dx):=\mu_t(dx)dt$-zero set is a $\mu^{(i)}(dt\,dx)=\mu^{(i)}_t(dx)dt$-zero set for both $i=1,2$.
Hence for $i=1,2$ by the Radon--Nikodym theorem there exist
$\mathcal B([s,T]\times H)$-measurable functions
$\rho_i\colon\, [s,T]\times H\to [0,\infty)$ such that
\begin{equation}
\label{e2.3}
\mu_t ^{(i)}(dx)dt:=\rho^{(i)}(t,x)\mu_t(dx)dt
\end{equation}
and it is easy to check from \eqref{e2.2} that $\rho^{(i)}\le 2$.
Furthermore, by Remark \ref{r1.1}(ii) for all $u\in D(L_0)$
$$
\int_s^T\int_HLu(t,x)\mu_t^{(1)}(dx)dt=\int_s^T\int_HLu(t,x)\mu_t^{(2)}(dx)dt,
$$
hence by \eqref{e2.3}
$$
\int_s^T\int_HLu(t,x)(\rho^{(1)}(t,x)-\rho^{(2)}(t,x)) \mu_t(dx)dt=0,\quad\forall\;  u\in D(L_0).
$$
Since $\mu$ satisfies \eqref{e2.1} and $\rho^{(1)}-\rho^{(2)}$ is bounded, it follows that $\rho^{(1)}-\rho^{(2)}=0,$ i.e. $\mu^{(1)}=\mu^{(2)}$. $\Box$\medskip

As we shall see in Sections 4 and 5 , sets as $\mathcal K$ arise very explicitly in the applications and are described by simple and natural integrability conditions.
\begin{Remark}
\label {r2.1}
\em Since we are in a parabolic situation Condition \eqref{e2.1} holds, if it holds with
$\lambda-L_0$ replacing $L_0$ for some $\lambda\in\R$.

\end{Remark}

\section{Regularity results for time dependent Ornstein--Uhlenbeck operators}

We need the following assumption on the coefficients  $A$ and $C$ in \eqref{e1.1}, \eqref{e1.1'}.
\begin{Hypothesis}
\label{h3.1}
\begin{enumerate}

\item[]

\item[{\rm(i)}] There is $\omega \in \R$ such that $ \langle Ax,x \rangle \le \omega
|x|^{2},\; \forall\;x\in D(A).$

\item[{\rm(ii)}]  $C\in L(H)$   is symmetric, nonnegative and such that the linear operator
$$
Q_t:=\int_0^t e^{sA}Ce^{sA^{*}}ds
$$
is of trace class for all $t>0$.

\item[{\rm(iii)}] One has  $e^{tA}(H)\subset Q_t^{1/2}(H)$ for all $t>0$
and there is    $\Lambda_t\in L(H)$ such that
$Q_t^{1/2}\Lambda_t=e^{tA}$ and
$$
\gamma_\lambda:=\int_0^{+\infty}e^{-\lambda t}\|\Lambda_t\|dt<+\infty,
$$
\end{enumerate}
where $\|\cdot\|$ denotes the operator norm in $L(H)$.
\end{Hypothesis}
By $R_t$ we denote the Ornstein--Uhlenbeck semigroup
$$
R_t\varphi(x):=\int_{H}^{}\varphi (e^{tA}x+y)N_{Q_t}(dy),\quad \varphi\in
C_{u,2}(H),
$$
where
$$
Q_tx:=\int_0^te^{sA}Ce^{sA^*}xds,\;\;\;x\in H, \quad t\ge 0,
$$
and $N_{Q_t}$ is the Gaussian measure in $H$ with mean $0$ and
covariance  operator $Q_t.$

We shall  consider $R_t$ acting in the Banach space $C_{u,2}(H)$, which consists of all   functions
  $\varphi\colon\, H\to \R$  such that the function
$x\mapsto\frac{\varphi(x)}{1+|x|^2}$ is uniformly continuous and bounded.
Let us define  the infinitesimal generator $U$ of $R_t$   through
its resolvent by setting, following \cite{Cerrai95},
$U:=\lambda-\widetilde{G_\lambda}^{-1},$
$D(U)=\widetilde{G_\lambda}(C_{u,2}(H))$, where
$$
\widetilde{G_\lambda}f(x)=\int_{0}^{+\infty }e^{-\lambda t}R_{t}f
(x)dt,\quad x\in H,\; \lambda>0,\;f\in C_{u,2}(H).
$$
It is easy to see that for any $h\in D(A^*)$ the  function  $\varphi_h(x)=e^{i\langle x,h
\rangle}$ belongs to the domain of $U$ in $C_{u,2}(H)$ and we have

\begin{equation}
\label{eA.1}
U\varphi_h = \frac{1}{2}\;\mbox{\rm Tr}\;[CD^{2}\varphi_h]+\langle x,A^*D\varphi_h
    \rangle.
\end{equation}
As a consequence of Hypothesis \ref{h3.1} one gets (see \cite[Lemma A.1]{BDPR08c})
\begin{Lemma}
\label{l3.2}
Let  Hypothesis \ref{h3.1} hold and let $\varphi\in D(U)$. Then there exists $c>0$ such that
$$
|D_x\varphi(x)|\le c(\|\varphi\|_{C_{u,2}(H)}+\|U\varphi\|_{C_{u,2}(H)})(1+|x|^2),\quad x\in H.
$$
\end{Lemma}
Now let us turn to the time-inhomogeneous case.  Let
$$
V_0u(t,x)=D_tu(t,x)+Uu(t,x),\quad u\in \mathcal
E_A([0,T]\times H).
$$
It is clear that $V_0u\in C([0,T];C_{u,2}(H))$ (note that $Uu(t,x)$ contains a
term growing as $|x|$).
Let us introduce an extension of the operator $V_0$.  For $\lambda\in \R$ set
$$
G_\lambda f(t,x)=\int_t^Te^{-\lambda(s-t)}R_{t-s}f(s,x)ds,\quad
f\in C([0,T];C_{u,2}(H)).
$$
It is easy to see that $G_\lambda$ satisfies the resolvent identity,   so that there exists
a unique linear closed operator $V$ in $C([0,T];C_{u,2}(H))$ such that
\begin{equation}
\label{eA.3}
G_\lambda =(\lambda-V)^{-1},\quad D(V)= G_\lambda(C([0,T];C_{u,2}(H))),\;\lambda\in \R.
\end{equation}
It is clear  that $V$ is an extension of $V_0$.

Finally, it is  easy to check that the semigroup $\mathcal R_\tau,\;\tau\ge 0,$ generated by the operator $V$
in the space $C_T([0,T];C_{u,2}(H)):=\{u\in C_T([0,T];C_{u,2}(H)):\;u(T,x)=0\}$ is given by

\begin{equation}
\label{eA.4}
\mathcal R_\tau f(t,x)=\left\{\begin{array}{l}
R_{\tau}f(t+\tau,\cdot)(x)\quad\mbox{\rm if}\;t+\tau\le T\\
0 \quad\mbox{\rm otherwise}.
\end{array}\right.
\end{equation}
Arguing as in \cite{Priola99} one can show  that $u\in D(V)$ and $Vu=f$ if and only if
\begin{equation}
\label{eA.5}
\left\{\begin{array}{l}
 \ds {\rm(i)}\;\;\;\;\lim\limits_{h\to 0}\frac1h\;
(\mathcal R_h u(t,x)-u(t,x))=f(t,x),\quad\forall\;(t,x)\in
[0,T]\times H,\\
\\
\ds {\rm(ii)}\;\;\; \sup\limits_{h\in (0,1],(t,x)\in
[0,T]\times H}\frac{(1+|x|^2)^{-1}}{h}\;|\mathcal R_h u(t,x)-u(t,x)|<+\infty.
\end{array}\right.
\end{equation}

We state now that $\mathcal E_{A}([0,T]\times H)$ is a   core for $V$.

The following results are generalization  of those in \cite{DPT01} and were proven in \cite[Proposition A.2, Corollary A3]{BDPR08c}.

\begin{Proposition}
\label{p3.3}
Let Hypothesis \ref{h3.1} hold and let $u\in D(V)$. Let $\nu$ be a finite
nonnegative Borel measure on~$[0,t]\times H$.
Then there exists a
sequence $(u_n)\subset {\mathcal E}_{A}([0,T]\times H)$
such that for some $c_1>0$ one has
$$
|u_n(t,x)|+|V_0u_n(t,x)|\le
c_1(1+|x|^2),\quad\forall\ (t,x)\in [0,T]\times H
$$
and $u_n\to u,$ $V_0u_n\to V_0u$ in measure~$\nu$.
 \end{Proposition}
 \begin{Corollary}
\label{c3.4}
Let Hypothesis \ref{h3.1} hold. Let $u\in D(V)$ and let $\nu$ be a finite nonnegative Borel measure
on $[0,T]\times H$.
Then there exists a sequence
$(u_n)\subset {\mathcal E}_{A}([0,T]\times H)$ such that for some
$c>0$ one has
$$
|u_n(t,x)|+|D_xu_n(t,x)|+|V_0u_n(t,x)|\le
c(1+|x|^2),\quad\forall\ (t,x)\in [0,T]\times H,
$$
and $u_n\to u,$ $D_xu_n\to D_xu,$ $V_0u_n\to Vu$
in measure~$\nu$.
 \end{Corollary}

We need the following
\begin{Hypothesis}
\label{h3.3}
$F\colon\, [0,T]\times H\to H$ is continuous
together with $D_xF(t,\cdot)\colon\, H\to L(H)$ for all $t\in [0,T].$ Moreover,
there exists $K>0$ such that
$$
|F(t,x)-F(t,y)|\le K|x-y|,\;\; x,y\in H,\; t\in [0,T].
$$
\end{Hypothesis}

Let $C_u(H,H)$ denote the set of all bounded uniformly continuous maps from $H$ to $H$ and $C^1_u(H)$ the set of all functions from $H$ to $\R$ which together with their first derivatives are bounded and uniformly continuous.

Then we have the following result  from \cite[Lemma 2.5]{BDPR08c}.

\begin{Proposition}
\label{p3.6}
Assume that Hypotheses \ref{h3.1} and \ref{h3.3} hold. Let $f\in C([0,T];C^1_u(H))$ and $\lambda\in\R$. Then there exists $u\in D(V)$ such that
\begin{enumerate}
\item[{\rm(i)}] $D_xu\in C([0,T];C_u(H,H))$,

\item[{\rm(ii)}] $\lambda u-Vu-\langle F,D_xu \rangle=f$,

\item[{\rm(iii)}] $\|u\|_\infty\le \frac1\lambda\;\|f\|_\infty,\quad \mbox{\rm if}\;\lambda>0.$
\end{enumerate}
\end{Proposition}

\section{The fully non-degenerate case}

In this section we shall consider the case where $C^{-1}\in L(H)$. Fix $\zeta\in\mathcal P(H)$, $s\in[0,T]$ and let
$\mathcal M_{s,\zeta}$ be defined as at the beginning of Section 2. The main result of this section is the following.
\begin{Theorem}
\label{t4.1}
Assume that Hypothesis \ref{h3.1} holds and that
\begin{equation}
\label{e4.1}
C^{-1}\in L(H).
\end{equation}
Define
$$
\mathcal K:=\left\{\mu\in \mathcal M_{s,\zeta}\colon\ \int_s^T\int_H(|x|^4+|F(t,x)|^2+|x|^4|F(t,x)|^2)\mu_t(dx)dt<\infty\right\}.
$$
(where again we set $|F(t,x)|=+\infty$ if $(t,x)\in [0,T]\times H\setminus D(F)$).
Then $\mathcal K$ contains at most one element.
\end{Theorem}
To prove Theorem \ref{t4.1}, by Theorem \ref{t2.1} we need to check that \eqref{e2.1} holds for all $\mu\in \mathcal K$. In fact we shall prove \eqref{e2.1} for an even  larger class of measures $\mathcal K^\lambda$ to be introduced below. We need some preparations. First for $\lambda\in [0,\infty)$ we introduce the set of measures
$\mathcal M_s^\lambda$ defined to be all finite nonnegative measures $\nu$ on $\mathcal B([0,T]\times H)$ satisfying \eqref{e1.2'} and \eqref{e1.5'} such that
\begin{equation}
\label{e4.2}
\int_s^T\int_HL_0u(t,x)\nu_t(dx)dt\le 2\lambda\int_s^T \int_Hu(t,x)\nu_t(dx)dt ,\quad\forall\;  u\in D(L_0),\;u\ge 0.
\end{equation}
\begin{Remark}
\label{r4.2}
\em By Remark \ref{r1.1}(ii) we have that $\mathcal M_{s,\zeta}\subset \mathcal M_s^\lambda$ for every $\zeta\in\mathcal P(H)$ and all $\lambda\ge 0$. Furthermore, we note that the set $\mathcal K$ defined in Theorem \ref{t4.1} is convex. For a large class of examples where $\mathcal K$ is nonempty  and thus consists of exactly one element we refer to Section 6.
\end{Remark}
\begin{Lemma}
\label{l4.3}
Let $\lambda\ge 0$ and $\nu\in \mathcal M_s^\lambda$ such that
\begin{equation}
\label{e4.3}
\int_s^T\int_H|F(t,x)|^2\nu_t(dx)dt<\infty.
\end{equation}
Then
\begin{equation}
\label{e4.4}
\begin{array}{l}
\ds \int_s^T\int_Hu(t,x) \,L_0u(t,x)\nu_t(dx)dt\le \lambda\int_s^T\int_Hu(t,x)^2 \nu_t(dx)dt\\
\\
\ds
-\frac12\;\int_s^T\int_H|\sqrt{C}\;D_xu(t,x)|^2\nu_t(dx)dt,\quad\forall\;  u\in D(L_0).
\end{array}
\end{equation}
In particular, $(L_0, D(L_0))$ is quasi-dissipative, hence closable in $L^2([0,T]\times H;\nu)$.
\end{Lemma}
{\bf Proof}. For all $u\in D(L_0)$ we have
\begin{equation}
\label{e4.5}
L_0 u^2=2u L_0u+|\sqrt{C}\;D_xu |^2,
\end{equation}
which implies \eqref{e4.4} by the definition of $\mathcal M_s^\lambda.$ The proof of the last part of the assertion is standard. $\Box$\medskip

For $\nu\in \mathcal M_s^\lambda$ satisfying \eqref{e4.3}, we denote the closure of $(L_0, D(L_0))$  on $L^2([0,T]\times H;\nu)$ by $(L^\nu, D(L^\nu))$. For $\lambda>0$ define
\begin{equation}
\label{e4.6}
\mathcal K^\lambda:=\left\{\nu\in \mathcal M^\lambda_{s}\colon\ \int_s^T\int_H(|x|^4+|F(t,x)|^2+|x|^4|F(t,x)|^2)\nu_t(dx)dt<\infty\right\}.
\end{equation}
Clearly, $\mathcal K \subset \mathcal K^\lambda$. Our aim is to prove that \eqref{e2.1} holds for all $\nu\in \mathcal K^\lambda$.
\begin{Remark}
\label{r4.4}
\em Once \eqref{e2.1} is proved for all $\nu\in \mathcal K^\lambda$, it follows that
$(L^\nu, D(L^\nu))$ is $m$-dissipative on $L^2([0,T]\times H;\nu)$, hence by the
Lumer--Phillips Theorem it generates a $C_0$-semigroup on $L^2([0,T]\times H;\nu)$.
However, we shall not use this fact below.
\end{Remark}

\begin{Lemma}
\label{l4.5}
Assume that Hypothesis \ref{h3.1} holds. Let $\lambda>0$ and $\nu\in \mathcal K^\lambda$ and let a map
 $F_0\colon\, [0,T]\times H\to H$ satisfy  Hypothesis \ref{h3.3}. Let $f\in C^1_u(H)$ and let $u_0$
 be as in Proposition~\ref{p3.3}, applied with $F_0$ replacing $F$, i.e. $u_0\in D(V)$, $\|u_0\|_\infty \le \frac1\lambda\;\|f\|_\infty$ and
$$
\lambda u_0-Vu_0- \langle F_0,D_x u_0  \rangle=f.
$$
Then:
\begin{enumerate}
\item[{\rm(i)}] $u_0\in D(L^\nu)$ and
\begin{equation}
\label{e4.7}
\lambda u_0-L^\nu u_0=f+ \langle F_0-F,D_x u_0  \rangle,
\end{equation}
as elements in $L^2([0,T]\times H;\nu)$.

\item[{\rm(ii)}] Suppose that \eqref{e4.1} holds. Then
\begin{multline*}
\int_s^T\int_H|D_xu(t,x)|^2\nu_t(dx)dt\\
\le \frac4\lambda\;\|C^{-1}\|\,\|f\|_\infty^2\Bigg((T-s)
 +\frac{\|C^{-1}\|}\lambda\;\int_s^T\int_H|F_0(t,x)-F(t,x)|^2\nu_t(dx)dt\Bigg).
\end{multline*}
\end{enumerate}
\end{Lemma}
{\bf Proof}. By Corollary \ref{c3.4} there exists $u_n\in D(L_0)$, $n\in \N$ such that
$$
u_n\to u_0,\quad D_xu_n\to D_xu_0,\quad V_0u _n\to Vu_0
$$
as $n\to \infty$ in $\nu$-measure and there exists $c\in (0,\infty)$ such that for all $(t,x)\in [0,T]\times H$
$$
|u_n(t,x)|+|V_0u_n(t,x)|+|D_xu_n(t,x)|\le c(1+|x|^2).
$$
Hence $L_0u_n\to Vu_0+\langle F,D_x u_0  \rangle$ as $n\to \infty$ in $\nu$-measure and
$$
|L_0u_n|\le c(1+|F(t,x)|)(1+|x|^2).
$$
Hence by assumption on $\nu$, Lebesgue's dominated convergence theorem implies that
$$
L_0u_n\to Vu_0+\langle F,D_x u_0  \rangle\quad\mbox{\rm as}\; n\to \infty\;\mbox{\rm in}\;L^2([0,T]\times H;\nu).
$$
Since $(L^\nu, D(L^\nu))$ is the closure of $(L^0, D(L^0))$, assertion (i) follows.

To prove (ii) we first note that by the above approximation and the assumptions on~$\nu$, \eqref{e4.4}
also holds  for~$u_0$. Hence multiplying \eqref{e4.7} by $u_0$ and integrating with respect to~$\nu$,
by the assumption on $\nu$ this implies
\begin{multline*}
\frac12\;\int_s^T\int_H|\sqrt{C}\;D_xu_0(t,x)|^2\nu_t(dx)dt
\le  \int_s^T\int_H|f(t,x)|\,|u_0(t,x)|\nu_t(dx)dt
\\
 + \int_s^T\int_H|F_0(t,x)-F(t,x)|\,|D_xu_0(t,x)|\,|u_0(t,x)|\nu_t(dx)dt.
\end{multline*}
Since $\|u_0\|_\infty\le \frac1\lambda\;\|f\|_\infty$, this implies assertion (ii). $\Box$\medskip

By Remark \ref{r4.2} and Theorem \ref{t2.1} the following result implies Theorem \ref{t4.1}.
\begin{Proposition}
\label{p4.6}
Let Hypothesis \ref{h3.1} and assumption \eqref{e4.1} hold. Let  $\lambda>0$ and $\nu\in \mathcal K^\lambda$. Then
\begin{equation}
\label{e4.8}
L_0(D(L_0))\;\mbox{\it is dense in }\;L^1([0,T]\times H,\nu),
\end{equation}
\end{Proposition}
{\bf Proof}. There exist $F_n\colon\, [0,T]\times H,\;n\in\N,$ satisfying Hypothesis \ref{h3.3}   such that
\begin{equation}
\label{e4.9}
\lim_{n\to\infty}\int_{[s,T]\times H}|F_n-F|^2d\nu=0.
\end{equation}
Let $f\in C^1_u(H)$ and $u_n$ as in Proposition \ref{p3.3}, applied with $F_n$ replacing $F$, i.e. $u_n\in D(V)$, $\|u_n\|_\infty\le \frac1\lambda\;\|f\|_\infty$ and
$$
\lambda u_n-V u_n-\langle F_n,D_xu_n   \rangle=f.
$$
Then by  Lemma \ref{l4.5}(i)
\begin{equation}
\label{e4.10}
\lambda u_n-L^\nu u_n=f+ \langle F_n-F,D_xu_n   \rangle
\end{equation}
and by Lemma \ref{l4.5}(ii) and \eqref{e4.9}
\begin{equation}
\label{e4.11}
\sup_{n\in\N}\int_s^T\int_H|D_xu_n(t,x)|^2\nu_t(dx)dt<\infty.
\end{equation}
\eqref{e4.9}-\eqref{e4.11} imply that $f$ is the closure of $(\lambda -L_0)(D(L_0))$ in
$L^1([0,T]\times H,\nu)$. Since $C_u^1(H)$ is dense in $L^1([0,T]\times H,\nu)$, \eqref{e4.8} now follows from Remark \ref{r2.1}. $\Box$

\section{Possibly degenerate cases}

In case $C$ is not invertible  and the noise is allowed to be very degenerate (including the deterministic $C=0$), more restrictive conditions on $F$ are needed to prove uniqueness of solutions to \eqref{e1.5}. Just for comparison with the results in the non degenerate case of the previous section we here recall the results from \cite{BDPR08c} which have no conditions on the noise and in particular include the case $C=0$.
\begin{Hypothesis}
\label{h5.1}
For each $t\in[0,T]$, $F(t,\cdot)$ is the minimal section of an $m$--dissipative graph
$$
\overline F(t,\cdot)\colon\, D(\overline F(t,\cdot))\subset H\to 2^H,\;t\in[0,T],
$$
i.e. for all $t\in[0,T]$, $D(\overline F(t,\cdot))\in\mathcal B(H)$ and there exists $K>0$ independent of $t$ such that
$$
 \langle u-v,x-y  \rangle\le K|x-y|^2,\quad\forall\;  x,y\in D(\overline F(t,\cdot)),\;u\in \overline F(t,x),\;v\in \overline F(t,y),
$$
and for every $\lambda>K$ one has
$$
\mbox{\rm Range}\;(\lambda-\overline F(t,\cdot)):=\bigcup_{x\in D(\overline F(t,\cdot))}
(\lambda x-\overline F(t,x))=H,$$
such that for all $t\in[0,T]$, $D(F(t,\cdot))=D(\overline F(t,\cdot))$ and for all
$x\in D(\overline F(t,\cdot))$, $F(t,x)\in \overline F(t,x)$ and $|F(t,x)|=\min_{y\in  \overline F(t,x)}|y|$. Furthermore, $0\in D(F(t,\cdot))$ and $F(t,0)=0$ for all $t\in[0,T]$.
\end{Hypothesis}

For $s\in[0,T]$, $\lambda>0$ let $\mathcal M_s^\lambda$ be as defined at the beginning of the previous section. Then the following result is proved in \cite[Theorem 3.3]{BDPR08c}.
 \begin{Theorem}
\label{t5.2}
Assume that Hypotheses \ref{h3.1} and \ref{h5.1} hold and let $s\in[0,T],\;\lambda>0$, $\nu\in \mathcal M_s^\lambda$  such that
$$
\int_s^T\int_H(|x|^2+|F(t,x)|+|x|^2|F(t,x)|)\nu_t(dx)dt<\infty.
$$
Then
\begin{equation}
\label{e5.1}
L_0(D(L_0))\;\mbox{\it is dense in }\;L^1([0,T]\times H,\nu),
\end{equation}
\end{Theorem}
\begin{Corollary}
\label{c4.3}
Assume Hypotheses \ref{h3.1} and \ref{h5.1}. Let $s\in[0,T],\;\zeta\in\mathcal P(H)$ and $M_{s,\zeta}$ be defined as at the beginning of Section 2. Define
\begin{equation}
\label{e5.2}
\mathcal K_1:=\left\{\mu\in \mathcal M_{s,\zeta}\colon\ \int_s^T\int_H(|x|^2+|F(t,x)|+|x|^2|F(t,x)|)\nu_t(dx)dt<\infty \right\}
\end{equation}
Then $\mathcal K_1$ contains at most one element.
\end{Corollary}
{\bf Proof}. Since $\mathcal K_1$  is convex the assertion follows immediately
from Theorems \ref{t5.2} and \ref{t2.1}. $\Box$

\begin{Remark}
\em  We note that once one assumes $F$ to satisfy Hypothesis \ref{h5.1} one can prove uniqueness in  $\mathcal K_1$ which is a larger set than $\mathcal K   $ in Theorem \ref{t4.1}, because of the weaker integrability condition. For  large classes of examples where $\mathcal K_1$ is non empty we refer to \cite{BDPR08b}, \cite{BDPR09} and \cite{BDPR10}.
\end{Remark}

\section{Applications}
Let $H=L^2(0,1):=L^2((0,1),d\xi)$  (with $|\cdot|:=|\cdot|_{L^2(0,1)}$) and let $A\colon\, D(A)\subset H\to H$ be defined by
$$
Ax(\xi)=\partial^2_\xi x(\xi),\;\xi\in (0,1),\quad D(A)=H^2(0,1)\cap H^1_0(0,1),
$$
where $\partial_\xi=\frac{d}{d\xi},$  $\partial^2_\xi=\frac{d^2}{d\xi^2}.$

We would like to mention here that what is done below generalizes to the case where $(0,1)$ is replaced by an open set $\mathcal O$ in $\R^d,\;d\ge 1.$ One has only to replace the operator $C$ below by $A^{-\delta}$ with properly chosen $\delta>0,$
depending on the dimension $d$.

Let $D(F):=[0,T]\times L^{2m}(0,1)$ and for $(t,\xi)\in D(F)$
$$
F(t,x)(\xi):=f(\xi,t,x(\xi))+h(\xi,t,x(\xi)),\quad\;\xi\in (0,1).
$$
Here $f,h\colon\, (0,1)\times[0,T]\times\R\to \R$ are   functions such that
 for every $\xi\in (0,1)$ the maps $f(\xi,\cdot,\cdot)$, $h(\xi,\cdot,\cdot)$ are
 continuous on $(0,T)\times\R$ and have the following properties:\medskip

\begin{enumerate}
\item[(f1)] (``polynomial bound''). There exist   $m\in\N$ and
a nonnegative function $c_1\in L^2(0,T)$  such that for all $t\in(0,T)$,
 $z\in\R,\;\xi\in (0,1)$ one has
$$
|f(\xi,t,z)|\le c_1(t)(1+|z|^m),
$$
also assuming without loss of generality that $m$ is odd.

\item[(f2)] (``quasi-dissipativity'').  There is a nonnegative function
$c_2\in L^1(0,T)$  such that for all $t\in[0,T],$ $z_1,z_2\in\R,\;\xi\in (0,1)$ one has
$$
(f(\xi,t,z_2)-f(\xi,t,z_1))(z_2-z_1)\le c_2(t)|z_2-z_1|^2.
$$

\item[(h1)] (``linear growth'').  There exists
a nonnegative function $c_3\in L^2(0,T)$ such that
for all $t\in[0,T],$ $z\in\R,\;\xi\in (0,1),$ one has
$$
|h(\xi,t,z)|\le c_3(t)(1+|z|).
$$
\end{enumerate}
Finally, let $C\in L(H)$ be symmetric, nonnegative and such that $C^{-1}\in L(H)$.

It is worth noting that it is not known whether
under these assumptions  the stochastic differential equation \eqref{e1.1} has a solution.

Set
\begin{equation}
\label{e6.1}
V_N(t,x):=\left\{\begin{array}{l}
2(c_1(t)+c_3(t)+1)(1+|x|^N_{L^{2N}(0,1)})\quad\mbox{\rm if}\;(t,x)\in [0,T]\times L^{2N}(0,1),\\
+\infty\quad\mbox{\rm otherwise}.
\end{array} \right.
\end{equation}
Observe, that by (f1) and (h1) one has
\begin{equation}
\label{e6.2}
|F(t,x)|\le V_m(t,x)<\infty\quad\forall\; (t,x)\in D(F).
\end{equation}
Let $N\ge m$. It was proved in \cite[Section 4]{BDPR10} that for every $\zeta\in \mathcal P(H)$  such that
\begin{equation}
\label{e6.3}
\int_H|x|^{2N}_{L^{2N}(0,1)}\zeta(dx)<+\infty
\end{equation}
(in particular for any Dirac measure with mass in $L^{2N}(0,1)$) there exists a solution $\mu(dt\,dx)=\mu_t(dx)dt$ to \eqref{e1.5}
satisfying \eqref{e1.2'}, \eqref{e1.5'} and in addition having the following properties
\begin{equation}
\label{e6.4}
\sup_{t\in[s,T]}\int_H|x|^{2}_{L^{2}(0,1)}\mu_t(dx)<\infty,
\end{equation}
\begin{equation}
\label{e6.5}
t\mapsto\int_Hu(t,x)\mu_t(dx)\;\mbox{\rm is continuous}\; \forall\;u\in D(L_0),
\end{equation}
\begin{equation}
\label{e6.6}
\begin{array}{l}
\ds\exists \;C>0\colon\ \int_s^T\int_H(V_N^2(r,x)+|(-A)^\delta x|^2_{L^{2}(0,1)})\mu_r(dx)dr\\
\\
\ds\le C\int_s^T\int_HV_N^2(r,x)\zeta(dx)dr<\infty ,\quad\forall\;\delta\in (\tfrac14, \tfrac12).
\end{array}
\end{equation}
Since by \eqref{e6.2} for $N:=m+2$ and some constant $C_1>0$ we have
$$
\begin{array}{l}
|x|^4_{L^2(0,1)}+|F(t,x)|^2_{L^2(0,1)}+|x|^4_{L^2(0,1)}\,|F(t,x)|^2_{L^2(0,1)}\\
\\
\le |x|^4_{L^{2N}(0,1)}+V^2_m(t,x)+|x|^4_{L^{2N}(0,1)}\,V^2_m(t,x)\\
\\
\le C_1V_N^2(t,x),
\end{array}
$$
it follows that, if
$$
\int_H |x|^{2(m+2)}_{L^{2(m+2)}(0,1)}\zeta(dx)<\infty,
$$
then the corresponding solution $\mu(dt\,dx)=\mu_t(dx)dt$ to \eqref{e1.5} is in the set $\mathcal K$ defined in Theorem \ref{t4.1} which in turn implies that it is the unique solution to   \eqref{e1.5} with $A,C,F$ as above such that
$$
\int_s^T\int_H(|x|^4+|F(t,x)|^2+|x|^4|F(t,x)|^2)\mu_t(dx)dt<\infty.
$$


\end{document}